\documentclass[12pt,reqno]{amsart}

\usepackage{amssymb,latexsym}

\usepackage{enumerate}

\usepackage[english]{babel}
\usepackage{amsmath}
\usepackage{graphicx}
\usepackage{amssymb}

\usepackage{amsthm,mathtools}

\usepackage{geometry}
\usepackage{tikz-cd}
\usepackage{mathrsfs}
\usepackage[colorinlistoftodos]{todonotes}
\usepackage{enumitem}
\usepackage{verbatim}
\usepackage[foot]{amsaddr}
\usepackage{xcolor}

\definecolor{myblue}{RGB}{0,90,180}
\makeatletter
\def\section{\@startsection{section}{1}%
  \z@{.7\linespacing\@plus\linespacing}{.5\linespacing}%
  {\normalfont\scshape\centering\color{myblue}}}
\def\subsection{\@startsection{subsection}{2}%
  \z@{.5\linespacing\@plus.7\linespacing}{-.5em}%
  {\normalfont\bfseries\color{myblue}}}
\def\subsubsection{\@startsection{subsubsection}{3}%
  \z@{.5\linespacing\@plus.7\linespacing}{-.5em}%
  {\normalfont\itshape\color{myblue}}}
\makeatother

\makeatletter

\@namedef{subjclassname@2010}{
	
	\textup{2020} Mathematics Subject Classification}

\makeatother
\newtheorem{thm}{Theorem}[section]
\newtheorem*{thm*}{Theorem}

\newtheorem{proposition}[thm]{Proposition}
\newtheorem{lemma}[thm]{Lemma}

\theoremstyle{definition}

\numberwithin{equation}{section}
\newcommand{\E}{\mathbb{E}}
\newcommand{\R}{\mathbb{R}}
\newcommand{\N}{\mathbb{N}}

\newcommand{\dd}{\,\mathrm{d}}

\newcommand{\p}{\mathfrak{p}}

\usepackage{hyperref}
\hypersetup{colorlinks=true,linkcolor=blue,anchorcolor=blue,citecolor=blue,
  pdftitle={High moments of zeta sums},
  pdfauthor={Zikang Dong, Weijia Wang, and Hao Zhang}}
\DeclareMathOperator\dif{d\!}

\newcommand{\newabstract}[1]{%
	\par\bigskip
	\csname otherlanguage*\endcsname{#1}%
	\csname captions#1\endcsname
	\item[\hskip\labelsep\scshape\abstractname.]
}

\begin{document}

	\baselineskip=17pt

	\title[High moments of zeta sums]{High moments of zeta sums}

	\author{Zikang Dong}
	\author{Weijia Wang}
	\author{Hao Zhang}
	\address[Zikang Dong]{School of Mathematical Sciences, Soochow University, Suzhou 215006, P. R. China}
	\address[Weijia Wang]{School of Mathematics, Shandong University, Jinan 250000, P. R. China}
	\address[Hao Zhang]{School of Mathematics, Hunan University, Changsha 410082, P. R. China}
	\email{zikangdong@gmail.com}
	\email{weijiawang@amss.ac.cn}
	\email{zhanghaomath@hnu.edu.cn}
	
	\date{}
	
	\begin{abstract}
For $2\le x\le T$, put $L=\min(x,T/x)$. We prove the unconditional
lower bound
\[
 \frac1T\int_T^{2T}\left|\sum_{n\le x}n^{-it}\right|^{2k}\dd t
 \gg_k x^k\bigl(\log(2L)\bigr)^{(k-1)^2}
\]
for each fixed real $k\ge2$. Under the Riemann hypothesis we obtain a
matching upper bound for every fixed real $k>2$, uniformly in $x$.
For positive integral $k$, the Steinhaus moment asymptotic of Harper,
Nikeghbali and Radziwi\l\l\ gives an explicit leading constant for
this time average when $x\to\infty$ and $x^k=o(T)$. Unconditionally, the fourth
moment has order $x^2\log(2L)$ throughout $2\le x\le T$. For each fixed integer $k\ge3$ and $0<\varepsilon<1/k$,
we also obtain the unconditional order
$x^k(\log(2T/x))^{(k-1)^2}$ for $T^{1-1/k+\varepsilon}\le x\le T$.
In particular, at $x=T/\log T$ the logarithmic factor in each of these
integer moments has order $(\log\log T)^{(k-1)^2}$.
\end{abstract}

	
	\maketitle
	
	\section{Introduction}
	
	%%%%%%%%%%%%%%%%%%%%%%%%%%%%%%%%%%%%%%%%%%%%%%%%%%%%%%%%%%%%%%%%%%%%%%%%%%%%%%
	Let $f(n)$ denote a Steinhaus random completely multiplicative function. Harper \cite{Harper19} obtained a breakthrough for low moments of partial sums of random multiplicative functions, both in the Steinhaus and Rademacher models. Character sums and zeta sums exhibit many features analogous to partial sums of Steinhaus random multiplicative functions. Harper \cite{Harper23} proved corresponding upper bounds for low moments of character and zeta sums. For lower bounds, we refer in particular to de la Bret\`eche, Munsch and Tenenbaum \cite{DMT}. Further background on zeta sums may be found in \cite{DWZ,Yang}.
    
    Harper \cite{Harper20} proved that for each fixed real $k\ge1$,
    \[
       \mathbb E\left|\sum_{n\le x}f(n)\right|^{2k}
       \asymp_k x^k(\log x)^{(k-1)^2}.
    \]
    For a large prime $q$ and $2\le x\le q/2$, put $L=\min(x,q/x)$.
    Szab\'o \cite[Theorem~1]{Szabo25} proved unconditionally, for real $k\ge2$, that
    \[
      \frac{1}{\varphi(q)}
      \sum_{\substack{\chi\bmod q\\ \chi\ne\chi_0}}
           \left|\sum_{n\le x}\chi(n)\right|^{2k}
      \gg_k x^k(\log L)^{(k-1)^2},
    \]
    extending earlier work of Munsch and Shparlinski \cite{MS16}.
    For $k>2$ the corresponding upper bound follows under the generalized
    Riemann hypothesis from \cite[Theorem~3]{Szabo24}, which treats
    primitive characters and builds on shifted-moment estimates
    (see also \cite{Munsch17}). For prime $q$ the primitive characters are
    precisely the non-principal ones, and in the above range
    $\min(\log x,\log(2q/x))\asymp\log L$.
    We refer to \cite{CZ98,GS01,MV79} for earlier results of moments of character sums.
	For high moments of zeta sums, Gao \cite{Gao24}, using shifted moment estimates for the Riemann zeta function proved by Curran \cite{Curran23}, established the following conditional upper bound.
	\begin{thm*}
		Let $k>2$ and $\varepsilon>0$ be fixed. Assume the Riemann hypothesis. Uniformly for $2\le x\le (1-\varepsilon)T$,
		\[
		\int_T^{2T}\bigg|\sum_{n\leq x}n^{-it}\bigg|^{2k}\,\dif t
		\ll_{k,\varepsilon} T x^k(\log T)^{(k-1)^2}.
		\]
	\end{thm*}
	For the rest of the paper, write
\[
  S_x(t)=\sum_{n\le x}n^{-it},\qquad
  M_k(x,T)=\frac1T\int_T^{2T}|S_x(t)|^{2k}\dd t,\qquad
  L=\min(x,T/x).
\]
Our first theorem is unconditional.
    
\begin{thm}\label{thm:main}
Fix a real number $k\ge2$. For all sufficiently large $T$, uniformly for
$2\le x\le T$, we have
\[
  \frac1T\int_T^{2T}\left|\sum_{n\le x}n^{-it}\right|^{2k}\,dt
  \gg_k x^k
  \begin{cases}
    (\log x)^{(k-1)^2},& 2\le x\le\sqrt T,\\[2pt]
    (\log(T/x))^{(k-1)^2},& \sqrt T<x\le T/2,\\[2pt]
    1,& T/2<x\le T.
  \end{cases}
\]
\end{thm}

Equivalently, Theorem~\ref{thm:main} gives
$M_k(x,T)\gg_k x^k(\log(2L))^{(k-1)^2}$.

\begin{thm}\label{thm:upper}
Assume the Riemann hypothesis. For each fixed real $k>2$, uniformly for
$2\le x\le T$ and all sufficiently large $T$,
\[
 M_k(x,T)\ll_k x^k\bigl(\log(2L)\bigr)^{(k-1)^2}.
\]
\end{thm}

Together, Theorems~\ref{thm:main} and~\ref{thm:upper} determine the
order of $M_k(x,T)$ under RH for every fixed real $k>2$. At
$x=T^\theta$ with $0<\theta<1$ fixed, Gao's upper bound already gives
this order, since $\log L=\min(\theta,1-\theta)\log T$. For lengths
such as $x=T/\log T$, Theorem~\ref{thm:upper} replaces the factor
$(\log T)^{(k-1)^2}$ in Gao's bound by a factor of order
$(\log\log T)^{(k-1)^2}$.

For a fixed positive integer $m$, Harper, Nikeghbali and Radziwi\l\l\
\cite[Theorem~3]{HNR} evaluated the Steinhaus moment
$\mathbb E|\sum_{n\le x}f(n)|^{2m}$ asymptotically, with an explicit
leading constant $C_m$. Heap and Lindqvist \cite{HL16} obtained this
independently. The Steinhaus expectation counts the diagonal solutions
$n_1\cdots n_m=n_{m+1}\cdots n_{2m}$. For the time average on $[T,2T]$,
the Montgomery--Vaughan mean-value theorem \cite{MV74} bounds the
off-diagonal contribution by $O(x^m/T)$ times this diagonal count
(see \eqref{eq:integer-transfer}). Hence
\begin{equation}\label{eq:integer-asymptotic}
 M_m(x,T)\sim C_m x^m(\log x)^{(m-1)^2}
 \qquad (x\to\infty,\ x^m=o(T)).
\end{equation}
In particular, $C_2=12/\pi^2$. We now give unconditional bounds for the
finite time average that also apply to long sums.

\begin{thm}\label{thm:integer}
Unconditionally, uniformly for $2\le x\le T$,
\[
 M_2(x,T)\asymp x^2\log(2L).
\]
For a fixed integer $m\ge3$ and fixed $0<\varepsilon<1/m$, the estimate
\[
 M_m(x,T)\asymp_{m,\varepsilon} x^m\bigl(\log(2L)\bigr)^{(m-1)^2}
\]
holds unconditionally for $2\le x\le T^{1/m}$ and for
$T^{1-1/m+\varepsilon}\le x\le T$.
\end{thm}

In the overlap $x\to\infty$, $x^m=o(T)$, formula
\eqref{eq:integer-asymptotic} identifies the leading constant in the order estimate of
Theorem~\ref{thm:integer}. At $x=T^{1/m}$ the relative error
bound in \eqref{eq:integer-transfer} is only $O(1)$.
The short-range upper bound in Theorem~\ref{thm:integer} still follows
from this estimate, and Theorem~\ref{thm:main} supplies the lower bound.

For long sums, Theorem~\ref{thm:integer} gives two-sided estimates on
$[T,2T]$ with a logarithm determined by the dual length $T/x$.
For the fourth moment these estimates run from $x=\sqrt T$ to $x=T$
and give the order $x^2\log(2T/x)$ throughout. For $m\ge3$ the unconditional
long range begins at $T^{1-1/m+\varepsilon}$. For example, the sixth
moment satisfies
\[
 M_3\!\left(\frac{T}{\log T},T\right)
 \asymp \left(\frac{T}{\log T}\right)^3(\log\log T)^4
\]
unconditionally. At this length $x^3/T\to\infty$, and the Steinhaus
sixth moment has order $x^3(\log T)^4$ by \cite[Theorem~3]{HNR}.
Thus the ratio of the time average to the random expectation has order
$(\log\log T/\log T)^4$ and tends to zero. The factor $\log\log T$ in
the time average comes from the lower limit $t/(2\pi x)\asymp T/x$
of the weighted dual sum in \eqref{eq:upper-flip}.
For $m\ge3$ the intervening range
$T^{1/m}<x<T^{1-1/m+\varepsilon}$ is covered by the upper bound of
Theorem~\ref{thm:upper} under RH.

To prove Theorem~\ref{thm:main} in Section~\ref{sec:lower}, we use the prime blocks of Szab\'o \cite[Sections~4 and~5]{Szabo25}
with cutoff $y=L^{1/C_0}$, where $C_0$ depends only on $k$.
In the time average the frequencies of a block polynomial are rational.
If their numerators and denominators are at most $H$, multiplying by the
original sum of length $x$ gives a relative mean-square error
$O(xH^2/T)$. This permits $H^2\le L^{1/2}$. In the random-model truncation,
conditioning on all primes outside the selected block eliminates a
factor $\log x/\log y$ from the error estimate
(Lemma~\ref{lem:weighted-block-moment}). This is needed when $T/x$
grows more slowly than every fixed power of $x$.

Section~\ref{sec:upper} proves Theorem~\ref{thm:upper}. On a dyadic block
of length $N$ we move to $\Re s=1/2+1/\log(2N)$ and use a mean-square
approximation to transfer Curran's shifted-moment bound to this line.
Integration over the shifts and the Carleson--Hunt theorem then give a
maximal estimate on the block. For $x>\sqrt T$ the approximate functional
equation starts the dual sum at $t/(2\pi x)$, producing the factor
$\log(2T/x)$. For Theorem~\ref{thm:integer}, a Dirichlet polynomial
mean-value estimate on the same blocks suffices.

\section{Proof of the lower bound}\label{sec:lower}

We first record a uniform second-moment lower bound
\begin{equation}\label{eq:baseline-long}
   \frac1T\int_T^{2T}|S_x(t)|^2\,dt\gg x
     \qquad (2\le x\le T).
\end{equation}
Fix a non-negative $w\in C_c^\infty((1,2))$ with $\int w=1$, and
let $\widehat w(u)=\int w(v)e^{-iuv}\,dv$. Choose a small absolute
$c>0$ and put $N=\max(1,\lfloor cx\rfloor)$. For every $m\le N$,
\[
 \sum_{\substack{n\le x\\n\ne m}}
       |\widehat w(T\log(n/m))|
    \ll_w (N/T)^2+x/T^2 .
\]
Indeed, for $n\le2N$ we have
$|\log(n/m)|\ge |n-m|/(2N)$ and may sum $|n-m|^{-2}$.
For $n>2N$ use $\log(n/m)\ge\log2$ and
$\widehat w(u)\ll_w(1+|u|)^{-2}$. On expanding the weighted scalar product,
the diagonal contributes $N$, while the off-diagonal is
$O_w(N((N/T)^2+x/T^2))$. Hence, with $c$ small and $T$ large,
\[
 \left|\frac1T\int_T^{2T}w(t/T)S_x(t)\overline{S_N(t)}\,dt\right|
     \ge N/2,\qquad
 \frac1T\int_T^{2T}w(t/T)|S_N(t)|^2\,dt\ll N .
\]
The second estimate follows from the same expansion with $x$ replaced
by $N$ in the sum over $n$. Cauchy--Schwarz and $N\gg x$ give
\eqref{eq:baseline-long}. Monotonicity of norms for $dt/T$ now gives
$T^{-1}\int_T^{2T}|S_x(t)|^{2k}\,dt\gg_k x^k$.

Set $L=\min(x,T/x)$ and assume henceforth that $L\ge L_0(k)$.
The case $L<L_0(k)$ follows from \eqref{eq:baseline-long}, since
$\log(2L)\le\log(2L_0(k))$. Choose $C_0=C_0(k)$ sufficiently large and
put $y=L^{1/C_0}$. In particular $x\ge L=y^{C_0}$ and $T/x\ge L$.
Put $y_0=1$, $y_M=y$ and $y_{m-1}=y_m^{1/20}$ for $2\le m\le M$.
Choose $M$ so that
\[
  (\log\log y)^2\le 20^{M-1}\le 20(\log\log y)^2.
\]
Increase $L_0(k)$ after fixing $C_0$ so that $M\ge2$. Define
\[
 J_M=\lfloor\sqrt{C_0}\rfloor,\qquad
 J_m=J_M+M-m\quad(2\le m<M),\qquad
 J_1=\lceil(\log\log y)^{3/2}\rceil .
\]
We may take $J_M\ge\exp(10^4k^2)$ and $J_1\ge J_M+1$. Since
$\log y_m=20^{m-M}\log y$, we have
\[
 \sum_{m=1}^M\frac{J_m}{20^{M-m}}
 \le \frac{J_1}{(\log\log y)^2}
       +\sum_{h\ge0}\frac{J_M+h}{20^h}
 \ll J_M+1 .
\]
Hence the logarithm of the product below, divided by $\log L$, is
$O(k(J_M+1)/C_0)$. Choosing $C_0$ large enough gives
\begin{equation}\label{eq:length-condition}
   \prod_{m=1}^M y_m^{10^4 kJ_m}\le L^{1/2}.
\end{equation}
All lower bounds imposed below on $C_0$ can be met simultaneously.
Once $C_0$ has been chosen, the remaining lower bounds on $y$ are secured by increasing $L_0(k)$.

For integers $\ell$ with
\[
    |\ell|\le \frac{\log y}{2},
\]
define
\begin{equation}\label{eq:Dml-t}
    D_{m,\ell}(t)
    :=
    \sum_{y_{m-1}<p\le y_m}
    \left(
        \frac{1}{p^{1/2+it+i\ell/\log y}}
        +
        \frac{1}{2p^{1+2it+2i\ell/\log y}}
    \right).
\end{equation}
Let
\[
    P_J(z):=\sum_{j=0}^J\frac{z^j}{j!},
\]
and define
\begin{equation}\label{eq:Rml-t}
    R_{m,\ell}(t)
    :=
    \left|
        P_{J_m}\bigl((k-1)\Re D_{m,\ell}(t)\bigr)
    \right|^2.
\end{equation}
We further set
\[
    A_\ell(t):=\prod_{m=1}^M R_{m,\ell}(t)
\]
and
\begin{equation}\label{eq:R-t}
    R(t):=
    \sum_{|\ell|\le (\log y)/2}A_\ell(t).
\end{equation}
Clearly,
\[
    R(t)\ge 0.
\]

\subsection{A mean-value lemma for rational frequencies}

The following lemma is the key point that allows us to transfer the
$t$-average to the Steinhaus random model.

\begin{lemma}\label{lem:rational-frequency}
Let
\[
    Q(t)=\sum_{\rho\in\mathcal R}c_\rho \rho^{-it},
\]
where the elements of $\mathcal R$ are distinct positive rational numbers.
Suppose that every $\rho\in\mathcal R$ can be written in reduced form as
\[
    \rho=\frac{a}{b},
    \qquad
    a\le A,\quad b\le B.
\]
Let $f$ be a Steinhaus random completely multiplicative function, and
define
\[
    Q(f)
    :=
    \sum_{\rho=a/b}
    c_\rho f(a)\overline{f(b)}.
\]
Then
\begin{equation}\label{eq:rational-transfer}
    \frac{1}{T}\int_T^{2T}|Q(t)|^2\,dt
    =
    \left(
        1+O\left(\frac{AB}{T}\right)
    \right)
    \mathbb E |Q(f)|^2.
\end{equation}
\end{lemma}

\begin{proof}
Let
\[
    \rho=\frac{a}{b},
    \qquad
    \rho'=\frac{c}{d}
\]
be two distinct elements of $\mathcal R$, both written in reduced form.
Then $ad\ne bc$, and hence
\[
    |ad-bc|\ge 1.
\]
Since $ad,bc\le AB$, we have
\[
    \left|\log \rho-\log \rho'\right|
    =
    \left|\log\frac{ad}{bc}\right|
    \gg \frac{1}{AB}.
\]
The logarithmic frequencies thus have separation
$\Delta\ge 1/(AB)$, by the mean-value theorem for the logarithm.
The Montgomery--Vaughan inequality \cite{MV74} gives
\[
    \int_T^{2T}|Q(t)|^2\,dt
    =
    T\sum_{\rho\in\mathcal R}|c_\rho|^2
    +
    O\left(
        \Delta^{-1}
        \sum_{\rho\in\mathcal R}|c_\rho|^2
    \right).
\]
Therefore
\begin{equation}\label{eq:rational-mv}
    \frac{1}{T}\int_T^{2T}|Q(t)|^2\,dt
    =
    \left(
        1+O\left(\frac{AB}{T}\right)
    \right)
    \sum_{\rho\in\mathcal R}|c_\rho|^2.
\end{equation}

On the other hand, the Steinhaus orthogonality relation gives
\[
    \mathbb E
    \left(
        f(a)\overline{f(b)}
        \overline{f(c)}f(d)
    \right)
    =
    \mathbf 1_{ad=bc}.
\]
Since the rational numbers $\rho$ in $\mathcal R$ are distinct and are
written in reduced form, it follows that
\[
    \mathbb E|Q(f)|^2
    =
    \sum_{\rho\in\mathcal R}|c_\rho|^2.
\]
Combining this with \eqref{eq:rational-mv} proves the lemma.
\end{proof}

\subsection{H\"older reduction}

By H\"older's inequality with exponents $k$ and $k/(k-1)$,
\begin{align}
    \frac{1}{T}\int_T^{2T}|S_x(t)|^2R(t)\,dt
    &\le
    \left(
        \frac{1}{T}\int_T^{2T}|S_x(t)|^{2k}\,dt
    \right)^{1/k}
    \nonumber\\
    &\quad\times
    \left(
        \frac{1}{T}\int_T^{2T}
        R(t)^{k/(k-1)}\,dt
    \right)^{(k-1)/k}.
    \label{eq:holder-main}
\end{align}
Consequently, it suffices to prove
\begin{equation}\label{eq:numerator-goal}
    \frac{1}{T}\int_T^{2T}|S_x(t)|^2R(t)\,dt
    \gg_k
    x(\log y)^{k^2-1},
\end{equation}
and
\begin{equation}\label{eq:denominator-goal}
    \frac{1}{T}\int_T^{2T}
    R(t)^{k/(k-1)}\,dt
    \ll_k
    (\log y)^{k^2+1}.
\end{equation}

\subsection{The numerator}

Let $f$ be a Steinhaus random completely multiplicative function. Define
\begin{equation}\label{eq:Dml-f}
    D_{m,\ell}(f)
    :=
    \sum_{y_{m-1}<p\le y_m}
    \left(
        \frac{f(p)}{p^{1/2+i\ell/\log y}}
        +
        \frac{f(p)^2}{2p^{1+2i\ell/\log y}}
    \right).
\end{equation}
Define $R_{m,\ell}(f)$, $A_\ell(f)$ and $R(f)$ by the same formulas as
above, replacing $D_{m,\ell}(t)$ by $D_{m,\ell}(f)$.

Put
\[
    F_{m,\ell}(t)
    :=
    P_{J_m}\bigl((k-1)\Re D_{m,\ell}(t)\bigr)
\]
and
\[
    F_\ell(t):=\prod_{m=1}^M F_{m,\ell}(t).
\]
Then
\[
    A_\ell(t)=|F_\ell(t)|^2.
\]

Every monomial occurring in $F_{m,\ell}(t)$ has rational frequency
$a/b$ with
\[
    a,b\le y_m^{2J_m}.
\]
Hence every rational frequency in $F_\ell(t)$ has numerator and
denominator at most
\[
    H:=\prod_{m=1}^M y_m^{2J_m}.
\]
It follows that
\[
    S_x(t)F_\ell(t)
\]
has rational frequencies whose reduced numerators and denominators are
bounded respectively by
\[
    xH
    \qquad\text{and}\qquad
    H.
\]
Thus Lemma~\ref{lem:rational-frequency} gives
\begin{align}
    \frac{1}{T}\int_T^{2T}
    |S_x(t)|^2A_\ell(t)\,dt
    &=
    \left(
        1+
        O\left(\frac{xH^2}{T}\right)
    \right)
    \mathbb E
    \left|
        \sum_{n\le x}f(n)
    \right|^2
    A_\ell(f).
    \label{eq:numerator-transfer-one}
\end{align}

By \eqref{eq:length-condition}, $H^2=\prod_m y_m^{4J_m}\le L^{1/2}$.
Since $L\le T/x$, this gives
\[
    \frac{xH^2}{T}\le \frac{L^{1/2}}{T/x}\le L^{-1/2}.
\]
We take $L_0(k)$ large enough that the relative error in
\eqref{eq:numerator-transfer-one} has absolute value at most $1/2$.
Summing \eqref{eq:numerator-transfer-one} over $\ell$, we obtain
\begin{equation}\label{eq:numerator-transfer}
    \frac{1}{T}\int_T^{2T}|S_x(t)|^2R(t)\,dt
    =
    (1+O(L^{-1/2}))
    \mathbb E
    \left|
        \sum_{n\le x}f(n)
    \right|^2R(f).
\end{equation}

We now use the corresponding random-model estimates. Put
\[
    E_\ell(f)
    :=
    \exp\left(
        2(k-1)
        \sum_{m=1}^M\Re D_{m,\ell}(f)
    \right).
\]
Proposition~\ref{prop:random-numerator-main} below gives
\begin{equation}\label{eq:random-main}
    \sum_{|\ell|\le(\log y)/2}
    \mathbb E
    \left|\sum_{n\le x}f(n)\right|^2 E_\ell(f)
    \gg_k x(\log y)^{k^2-1}.
\end{equation}
For each block put
\[
    \operatorname{Err}_{m,\ell}(f)
    :=e^{2(k-1)\Re D_{m,\ell}(f)}-R_{m,\ell}(f).
\]
Proposition~\ref{prop:random-truncation} below shows that
\begin{align}
&\sum_{|\ell|\le(\log y)/2}
\mathbb E\left|\sum_{n\le x}f(n)\right|^2
\exp\left(2(k-1)\sum_{m'\ne m}\Re D_{m',\ell}(f)\right)
|\operatorname{Err}_{m,\ell}(f)|
\nonumber\\
&\qquad\ll_k e^{-J_m}
 x(\log y)^{k^2-1}.
\label{eq:random-truncation}
\end{align}
Now
\[
R(f)=\sum_{|\ell|\le(\log y)/2}E_\ell(f)
\prod_{m=1}^M\left(1-
\frac{\operatorname{Err}_{m,\ell}(f)}
{e^{2(k-1)\Re D_{m,\ell}(f)}}\right).
\]
Each factor in this product is non-negative. Hence, using
$\prod_m(1+u_m)\ge1-\sum_m|u_m|$ whenever $u_m\ge-1$, we obtain
\begin{align*}
\mathbb E\left|\sum_{n\le x}f(n)\right|^2R(f)
&\ge \sum_{|\ell|\le(\log y)/2}
\mathbb E\left|\sum_{n\le x}f(n)\right|^2E_\ell(f)\\
&\quad-\sum_{m=1}^M\sum_{|\ell|\le(\log y)/2}
\mathbb E\left|\sum_{n\le x}f(n)\right|^2
 e^{2(k-1)\sum_{m'\ne m}\Re D_{m',\ell}(f)}
 |\operatorname{Err}_{m,\ell}(f)|.
\end{align*}
Since $J_1\ge J_M+1$, summing the geometric progression gives
$\sum_m e^{-J_m}\le 2e^{-J_M}$. The ratio of the total error bound
to $x(\log y)^{k^2-1}$ is therefore $O_k(e^{-\lfloor\sqrt{C_0}\rfloor})$.
The lower-bound constant in \eqref{eq:random-main} is independent of $C_0$
once $C_0$ exceeds a threshold depending only on $k$.
Enlarging $C_0=C_0(k)$ gives
\begin{equation}\label{eq:random-numerator}
    \mathbb E\left|\sum_{n\le x}f(n)\right|^2R(f)
    \gg_k x(\log y)^{k^2-1}.
\end{equation}
Combining \eqref{eq:numerator-transfer} and
\eqref{eq:random-numerator}, we obtain
\[
    \frac{1}{T}\int_T^{2T}|S_x(t)|^2R(t)\,dt
    \gg_k
    x(\log y)^{k^2-1},
\]
which proves \eqref{eq:numerator-goal}.

\subsection{The denominator}

Put
\[
    q:=\frac{k}{k-1}
    =
    1+\frac{1}{k-1}.
\]
Since $k\ge2$, the exponent $1/(k-1)$ lies in $(0,1]$.
Subadditivity (with equality when $k=2$) gives
\[
    \left(
        \sum_{\ell'}A_{\ell'}(t)
    \right)^{1/(k-1)}
    \le
    \sum_{\ell'}
    A_{\ell'}(t)^{1/(k-1)}.
\]
Thus
\begin{equation}\label{eq:Rq-double-sum}
    R(t)^q
    \le
    \sum_{\ell,\ell'}
    \prod_{m=1}^M
    R_{m,\ell}(t)
    R_{m,\ell'}(t)^{1/(k-1)}.
\end{equation}

For each $m$, put
\[
    I_0^{(m)}
    :=
    \left[0,\frac{J_m}{100k}\right)
\]
and, for $n\ge1$,
\[
    I_n^{(m)}
    :=
    \left[
        \frac{2^{n-1}J_m}{100k},
        \frac{2^nJ_m}{100k}
    \right).
\]
For a fixed tuple
\[
    \mathbf n=(n_1,\ldots,n_M),
\]
put
\[
    W_m:=\inf I_{n_m}^{(m)}.
\]
We decompose according to the conditions
\[
    |\Re D_{m,\ell'}(t)|
    \in I_{n_m}^{(m)}.
\]

Let
\[
    a_m:=2\lceil200kJ_m\rceil.
\]
For the fixed dyadic cell define
\begin{equation}\label{eq:def-U}
U_{m,\ell'}^{(n_m)}(t):=
\begin{cases}
\displaystyle
\left|P_{J_m}(\Re D_{m,\ell'}(t))\right|^2, & n_m=0,\\[1.5ex]
\displaystyle
e^{4W_m}\left|\frac{D_{m,\ell'}(t)}{W_m}\right|^{a_m},
& \dfrac{J_m}{100k}\le W_m\le100kJ_m,\\[2ex]
\displaystyle
\left(\frac{2(k-1)^{J_m}(2W_m)^{J_m}}{J_m!}\right)^{2/(k-1)}
\left|\frac{D_{m,\ell'}(t)}{W_m}\right|^{a_m},
& W_m>100kJ_m.
\end{cases}
\end{equation}
On the good cell, $|v|\le J_m/(100k)$ for $v=\Re D_{m,\ell'}(t)$.
Taylor's remainder and $(J+1)!\ge((J+1)/e)^{J+1}$ give
$P_{J_m}(z)=e^z(1+O(e^{-J_m}))$ for both $z=v$ and $z=(k-1)v$.
On a middle cell, $P_{J_m}((k-1)v)$ has absolute value at most
$e^{(k-1)|v|}$, and $|D_{m,\ell'}|\ge W_m$. On a large cell the
successive terms of $\sum_{j\le J_m}((k-1)|v|)^j/j!$, read backwards
from $j=J_m$, decrease at least geometrically with ratio $1/2$.
Since $|v|\le 2W_m$, these estimates give
\begin{equation}\label{eq:U-majorant}
    R_{m,\ell'}(t)^{1/(k-1)}
    \le (1+O(e^{-J_m}))U_{m,\ell'}^{(n_m)}(t).
\end{equation}
The error is uniform in $m$ and in the cell for fixed $k$. In particular,
\[
    \prod_{m=1}^M(1+O(e^{-J_m}))\ll1,
\]
because $\sum_m e^{-J_m}\ll e^{-J_M}$. Each majorant in \eqref{eq:def-U} is the squared modulus of a finite rational-frequency polynomial; the degree of an underlying polynomial is $O(kJ_m)$.

Consequently, for fixed $\ell,\ell'$ and $\mathbf n$,
\[
    \prod_{m=1}^M
    R_{m,\ell}(t)
    U_{m,\ell'}^{(n_m)}(t)
\]
is the squared modulus of a rational-frequency polynomial. Each block of
its square root contributes numerator and denominator at most
$y_m^{2J_m+\max(2J_m,a_m)}\le y_m^{1000kJ_m}$, since the terms of
$D_{m,\ell'}$ involve $p$ or $p^2$ and $a_m$ is even. Put
\[
    G=\prod_{m=1}^M y_m^{1000kJ_m}.
\]
Both heights are bounded by $G$ independently of the cell and the shifts.
By \eqref{eq:length-condition}, $G^2\le L^{1/2}$ and
$G^2/T\le T^{-3/4}$ since $L\le\sqrt T$.
Lemma~\ref{lem:rational-frequency} therefore yields
\begin{align}
    &\frac{1}{T}
    \int_T^{2T}
    \prod_{m=1}^M
    R_{m,\ell}(t)
    U_{m,\ell'}^{(n_m)}(t)\,dt
    \nonumber\\
    &\qquad=
    (1+O(T^{-3/4}))
    \mathbb E
    \prod_{m=1}^M
    R_{m,\ell}(f)
    U_{m,\ell'}^{(n_m)}(f).
    \label{eq:denominator-transfer}
\end{align}

Proposition~\ref{prop:mixed-random-complete} below gives, uniformly in $\ell,\ell'$ and the dyadic cell,
\begin{equation}\label{eq:mixed-random}
    \mathbb E\prod_{m=1}^M R_{m,\ell}(f)U_{m,\ell'}^{(n_m)}(f)
    \ll_k \frac{(\log y)^{k^2}}{1+|\ell-\ell'|^{2(k-1)}}
    \prod_{m=1}^M(1+W_m)^{-2}.
\end{equation}

Apply \eqref{eq:U-majorant} on each cell in \eqref{eq:Rq-double-sum}.
The product of the majorants is non-negative, so its integral over a
cell is at most its integral over $[T,2T]$. Summing
\eqref{eq:denominator-transfer} and \eqref{eq:mixed-random} now gives
\begin{align}
    \frac{1}{T}\int_T^{2T}R(t)^q\,dt
    &\ll_k
    (\log y)^{k^2}
    \sum_{\ell,\ell'}
    \frac{1}
    {1+|\ell-\ell'|^{2(k-1)}}
    \nonumber\\
    &\qquad\times
    \sum_{\mathbf n}
    \prod_{m=1}^M
    (1+W_m)^{-2}.
    \label{eq:denominator-sums}
\end{align}

For every $m$, the dyadic definition gives
\[
\sum_{n_m\ge0}\left(1+\inf I_{n_m}^{(m)}\right)^{-2}
\le 1+2\left(\frac{100k}{J_m}\right)^2.
\]
Consequently
\begin{align*}
\sum_{\mathbf n}\prod_{m=1}^M(1+W_m)^{-2}
&=\prod_{m=1}^M\sum_{n_m\ge0}\left(1+\inf I_{n_m}^{(m)}\right)^{-2}\\
&\le\exp\left(10^5k^2\sum_{m=1}^M J_m^{-2}\right)\ll_k1,
\end{align*}
because $J_m=J_M+M-m$ for $m\ge2$ and $J_1\to\infty$.
Moreover, since $2(k-1)\ge2$,
\[
    \sum_{h\in\mathbb Z}
    \frac{1}{1+|h|^{2(k-1)}}
    \ll_k 1.
\]
There are $O(\log y)$ possible values of $\ell$, and consequently
\[
    \sum_{|\ell|,|\ell'|\le(\log y)/2}
    \frac{1}{1+|\ell-\ell'|^{2(k-1)}}
    \ll_k \log y.
\]
It follows from \eqref{eq:denominator-sums} that
\[
    \frac{1}{T}\int_T^{2T}
    R(t)^{k/(k-1)}\,dt
    \ll_k
    (\log y)^{k^2+1}.
\]
This proves \eqref{eq:denominator-goal}.

\subsection{Completion of the proof}

\begin{proof}[Proof of Theorem~\ref{thm:main}]
By \eqref{eq:holder-main},
\eqref{eq:numerator-goal} and \eqref{eq:denominator-goal},
\begin{align*}
    \frac{1}{T}\int_T^{2T}|S_x(t)|^{2k}\,dt
    &\gg_k
    \frac{
        \left(
            x(\log y)^{k^2-1}
        \right)^k
    }{
        \left(
            (\log y)^{k^2+1}
        \right)^{k-1}
    }\\
    &=
    x^k
    (\log y)^{
        k(k^2-1)-(k-1)(k^2+1)
    }.
\end{align*}
Since
\[
    k(k^2-1)-(k-1)(k^2+1)
    =(k-1)^2,
\]
we obtain
\[
    \frac{1}{T}\int_T^{2T}|S_x(t)|^{2k}\,dt
    \gg_k
    x^k(\log y)^{(k-1)^2}.
\]
For $L\ge L_0(k)$ we have $\log y=(\log L)/C_0(k)\asymp_k\log(2L)$.
Combining this with \eqref{eq:baseline-long} for bounded $L$ yields
\[
    \frac1T\int_T^{2T}|S_x(t)|^{2k}\,dt
    \gg_k x^k\bigl(\log(2\min(x,T/x))\bigr)^{(k-1)^2}
    \qquad (2\le x\le T).
\]
This completes the proof.
\end{proof}

\subsection{Random-model estimates}

The untruncated numerator and its error correspond to
\cite[Propositions~4.1 and~4.2]{Szabo25}. For the denominator we use the
good- and bad-block decomposition in \cite[Section~5]{Szabo25}. We give the
estimates with the truncation orders chosen in \eqref{eq:length-condition}.  Throughout this subsection, $f$ denotes a
Steinhaus random completely multiplicative function.

For convenience, write
\[
    F_y(s)
    :=
    \prod_{p\le y}
    \left(1-\frac{f(p)}{p^s}\right)^{-1}.
\]

\subsubsection{A two-point Euler-product estimate}

We first record a standard calculation for correlated Euler products.

\begin{lemma}\label{lem:euler-correlation}
Let $\alpha,\gamma\ge0$ be fixed, and let
$\sigma_1,\sigma_2\ge0$ and $t_1,t_2\in\mathbb R$. Then
\begin{align}
&\mathbb E
\prod_{z<p\le y}
\left|
1-\frac{f(p)}{p^{1/2+\sigma_1+it_1}}
\right|^{-2\alpha}
\left|
1-\frac{f(p)}{p^{1/2+\sigma_2+it_2}}
\right|^{-2\gamma}
\nonumber\\
&\quad=
\exp\left(
\sum_{z<p\le y}
\left\{
\frac{\alpha^2}{p^{1+2\sigma_1}}
+
\frac{\gamma^2}{p^{1+2\sigma_2}}
+
\frac{
2\alpha\gamma
\cos((t_1-t_2)\log p)
}{
p^{1+\sigma_1+\sigma_2}
}
+
O_{\alpha,\gamma}(p^{-3/2})
\right\}
\right),
\label{eq:euler-correlation}
\end{align}
provided $z$ is sufficiently large in terms of $\alpha$ and $\gamma$.
Taking $z=z(\alpha,\gamma)$ fixed, the omitted primes contribute a
factor bounded above and below by positive constants depending only on
$\alpha$ and $\gamma$, since $|f(p)p^{-1/2-\sigma-it}|\le p^{-1/2}$.
\end{lemma}

\begin{proof}
By independence of the random variables $f(p)$, the expectation
factorizes prime by prime. Put
\[
    u_p=\frac{f(p)}{p^{1/2+\sigma_1+it_1}},
    \qquad
    v_p=\frac{f(p)}{p^{1/2+\sigma_2+it_2}}.
\]
Using
\[
    -\log(1-z)=z+\frac{z^2}{2}+O(|z|^3)
\]
uniformly for $|z|\le1/2$, the local factor equals
\[
\exp\Big(2\alpha\Re u_p+\alpha\Re(u_p^2)
          +2\gamma\Re v_p+\gamma\Re(v_p^2)
          +O_{\alpha,\gamma}(p^{-3/2})\Big).
\]
Put $X_p=2\alpha\Re u_p+2\gamma\Re v_p$. Since $X_p=O_{\alpha,\gamma}(p^{-1/2})$, expansion to second order gives
\[
1+X_p+\alpha\Re(u_p^2)+\gamma\Re(v_p^2)+\frac{X_p^2}{2}
+O_{\alpha,\gamma}(p^{-3/2}).
\]
Averaging over the unit circle, the terms linear in $f(p)$ and $f(p)^2$ vanish. Moreover,
\[
    \mathbb E(2\alpha\Re u_p)^2/2
    =
    \frac{\alpha^2}{p^{1+2\sigma_1}},
\]
\[
    \mathbb E(2\gamma\Re v_p)^2/2
    =
    \frac{\gamma^2}{p^{1+2\sigma_2}},
\]
and
\[
    \mathbb E
    (2\alpha\Re u_p)(2\gamma\Re v_p)
    =
    \frac{
        2\alpha\gamma
        \cos((t_1-t_2)\log p)
    }{
        p^{1+\sigma_1+\sigma_2}
    }.
\]
All terms of total degree at least three contribute
$O_{\alpha,\gamma}(p^{-3/2})$. Hence the local factor equals
\[
    1+
    \frac{\alpha^2}{p^{1+2\sigma_1}}
    +
    \frac{\gamma^2}{p^{1+2\sigma_2}}
    +
    \frac{
        2\alpha\gamma
        \cos((t_1-t_2)\log p)
    }{
        p^{1+\sigma_1+\sigma_2}
    }
    +
    O_{\alpha,\gamma}(p^{-3/2}).
\]
Multiplying over the primes and using
$\log(1+u)=u+O(u^2)$ proves the lemma.
\end{proof}

We shall use the following standard shifted prime-harmonic estimate. Uniformly for $t\in\mathbb R$ and $y\ge2$,
\begin{equation}\label{eq:prime-harmonic-shift}
\sum_{p\le y}\frac{\cos(t\log p)}p\le
\begin{cases}
\log\log y+O(1), & |t|\le1/\log y,\\
\log(1/|t|)+O(1), & 1/\log y\le |t|\le10,\\
\log\log|t|+O(1), & |t|\ge10.
\end{cases}
\end{equation}
In addition, if $|t|\le1/\log y$, then
\begin{equation}\label{eq:prime-harmonic-lower}
\sum_{p\le y}\frac{\cos(t\log p)}{p^{1+\eta}}
\ge \log\log y-\log(2+\eta\log y)+O(1)
\end{equation}
uniformly for $0\le\eta\le1$. For \eqref{eq:prime-harmonic-shift}, see \cite[Lemma~3]{Szabo25}.
For the lower bound, the contribution of $1-\cos(t\log p)$ is
$O(t^2\sum_{p\le y}(\log p)^2/p)=O(1)$. Restrict the remaining sum
to $p\le \exp(\min(\log y,1/\eta))$, where replacing $p^{-\eta}$ by $1$
costs $O(\eta\sum (\log p)/p)=O(1)$. Mertens' theorem now gives
\eqref{eq:prime-harmonic-lower}, with $1/\eta=+\infty$ for $\eta=0$.

\subsubsection{The random-model numerator}

We prove
\begin{proposition}\label{prop:random-numerator-main}
For fixed $k\ge2$ and all $C_0$ sufficiently large in terms of $k$,
\begin{align}
&\sum_{|\ell|\le(\log y)/2}
\mathbb E
\left|
    \sum_{n\le x}f(n)
\right|^2
\exp\left(
    2(k-1)
    \Re
    \sum_{p\le y}
    \left(
        \frac{f(p)}{p^{1/2+i\ell/\log y}}
        +
        \frac{f(p)^2}{2p^{1+2i\ell/\log y}}
    \right)
\right)
\nonumber\\
&\hspace{5cm}
\gg_k
x(\log y)^{k^2-1}.
\label{eq:random-main-complete}
\end{align}
\end{proposition}

\begin{proof}
For every $|\ell|\le(\log y)/2$,
\[
\sum_{p\le y}
\left(
    \frac{f(p)}{p^{1/2+i\ell/\log y}}
    +
    \frac{f(p)^2}{2p^{1+2i\ell/\log y}}
\right)
=
\log F_y\left(
    \frac12+\frac{i\ell}{\log y}
\right)
+
O(1),
\]
uniformly in $f$ and $\ell$, since
\[
    \sum_p p^{-3/2}<\infty.
\]
Consequently,
\begin{align}
&\exp\left(
    2(k-1)\Re
    \sum_{p\le y}
    \left(
        \frac{f(p)}{p^{1/2+i\ell/\log y}}
        +
        \frac{f(p)^2}{2p^{1+2i\ell/\log y}}
    \right)
\right)
\nonumber\\
&\qquad\asymp_k
\left|
F_y\left(
    \frac12+\frac{i\ell}{\log y}
\right)
\right|^{2(k-1)}.
\label{eq:exp-to-Euler}
\end{align}

Let $\mathbb E^{(y)}$ denote conditional expectation given the variables
$(f(p))_{p\le y}$. We first need a lower bound for
\[
    \mathbb E^{(y)}
    \left|
        \sum_{n\le x}f(n)
    \right|^2.
\]

Write $P^+(n)$ and $P^-(n)$ for the largest and smallest prime factors
of $n$, with $P^+(1)=1$ and $P^-(1)=\infty$. Every positive integer $n$
has a unique factorization
\[
    n=ab,
    \qquad
    P^+(a)\le y,
    \qquad
    P^-(b)>y.
\]
Hence
\begin{align*}
\mathbb E^{(y)}
\left|
    \sum_{n\le x}f(n)
\right|^2
&=
\sum_{\substack{b\le x\\P^-(b)>y}}
\left|
    \sum_{\substack{a\le x/b\\P^+(a)\le y}}
    f(a)
\right|^2.
\end{align*}
Restricting $b$ to a suitable subrange and putting
$r=\lfloor x/b\rfloor$, we obtain
\begin{align}
\mathbb E^{(y)}
\left|
    \sum_{n\le x}f(n)
\right|^2
&\ge
\sum_{r\le x^{1/4}}
\left|
    \sum_{\substack{a\le r\\P^+(a)\le y}}f(a)
\right|^2
\sum_{\substack{x/(r+1)<b\le x/r\\P^-(b)>y}}1.
\label{eq:conditioned-sieve-1}
\end{align}

Write $L_r=x/(r(r+1))\gg x^{1/2}$ for $r\le x^{1/4}$.
The number of multiples of a squarefree $d$ in the interval is
$L_r/d+O(1)$. Apply the one-dimensional lower-bound sieve
\cite{FI} to the primes $p\le y$ with level $D=y^3$.
Here $\log D/\log y=3$, so the lower sieve function has a fixed positive value.
The remainder is $O(D)=o(L_r/\log y)$ since $y\le x^{1/C_0}$ and $C_0\ge 10$.
Mertens' theorem thus gives, with an absolute lower-bound constant,
\[
    \sum_{\substack{x/(r+1)<b\le x/r\\P^-(b)>y}}1
      \gg \frac{L_r}{\log y}\gg\frac{x}{r^2\log y}
      \qquad (r\le x^{1/4}).
\]
Thus
\begin{align}
\mathbb E^{(y)}
\left|
    \sum_{n\le x}f(n)
\right|^2
&\gg
\frac{x}{\log y}
\sum_{r\le x^{1/4}}
\frac1{r^2}
\left|
    \sum_{\substack{a\le r\\P^+(a)\le y}}f(a)
\right|^2.
\label{eq:conditioned-sieve-2}
\end{align}
Partial summation yields
\begin{align}
\mathbb E^{(y)}
\left|
    \sum_{n\le x}f(n)
\right|^2
&\gg
\frac{x}{\log y}
\int_1^{x^{1/4}}
\left|
    \sum_{\substack{a\le u\\P^+(a)\le y}}f(a)
\right|^2
\frac{du}{u^2}.
\label{eq:conditioned-integral}
\end{align}

Let $0<\beta<1/10$. By the usual Rankin trick,
\begin{align}
\int_1^{x^{1/4}}
\left|
    \sum_{\substack{a\le u\\P^+(a)\le y}}f(a)
\right|^2
\frac{du}{u^2}
&\ge
\int_1^\infty
\left|
    \sum_{\substack{a\le u\\P^+(a)\le y}}f(a)
\right|^2
\frac{du}{u^{2+2\beta}}
\nonumber\\
&\quad-
x^{-\beta/4}
\int_1^\infty
\left|
    \sum_{\substack{a\le u\\P^+(a)\le y}}f(a)
\right|^2
\frac{du}{u^{2+\beta}}.
\label{eq:rankin-split}
\end{align}
By the standard Parseval formula for Dirichlet series,
\begin{align}
\int_1^\infty
\left|
    \sum_{\substack{a\le u\\P^+(a)\le y}}f(a)
\right|^2
\frac{du}{u^{2+2\beta}}
&\gg
\int_{-1/2}^{1/2}
|F_y(1/2+\beta+it)|^2\,dt,
\end{align}
whereas
\begin{align}
\int_1^\infty
\left|
    \sum_{\substack{a\le u\\P^+(a)\le y}}f(a)
\right|^2
\frac{du}{u^{2+\beta}}
&\ll
\int_{-\infty}^{\infty}
\frac{
|F_y(1/2+\beta/2+it)|^2
}{
1/4+t^2
}\,dt.
\end{align}
Therefore
\begin{align}
\mathbb E^{(y)}
\left|
    \sum_{n\le x}f(n)
\right|^2
&\gg
\frac{x}{\log y}
\Bigg\{
\int_{-1/2}^{1/2}
|F_y(1/2+\beta+it)|^2\,dt
\nonumber\\
&\qquad\qquad-
x^{-\beta/4}
\int_{-\infty}^{\infty}
\frac{
|F_y(1/2+\beta/2+it)|^2
}{
1/4+t^2
}\,dt
\Bigg\}.
\label{eq:conditioned-final}
\end{align}

We multiply \eqref{eq:conditioned-final} by
\[
    \left|
    F_y\left(
        \frac12+\frac{i\ell}{\log y}
    \right)
    \right|^{2(k-1)}
\]
and take expectation.

Take $\beta=1/\log y$, which lies in $(0,1/10)$ once $L_0(k)$
is sufficiently large.

Consider first
\begin{align}
\mathcal I_\ell
:=
\int_{-1/2}^{1/2}
\mathbb E
|F_y(1/2+\beta+it)|^2
\left|
F_y\left(
    \frac12+\frac{i\ell}{\log y}
\right)
\right|^{2(k-1)}
dt.
\end{align}
By Lemma~\ref{lem:euler-correlation},
\begin{align}
&\mathbb E
|F_y(1/2+\beta+it)|^2
\left|
F_y\left(
    \frac12+\frac{i\ell}{\log y}
\right)
\right|^{2(k-1)}
\nonumber\\
&\quad\asymp_k
\exp\left(
\sum_{p\le y}
\frac1{p^{1+2\beta}}
+
(k-1)^2\sum_{p\le y}\frac1p
+
2(k-1)
\sum_{p\le y}
\frac{
\cos((t-\ell/\log y)\log p)
}{
p^{1+\beta}
}
\right).
\label{eq:first-correlation}
\end{align}

If
\[
    \left|
        t-\frac{\ell}{\log y}
    \right|
    \le \frac1{\log y},
\]
then
\[
    \sum_{p\le y}\frac1{p^{1+2\beta}}
    =
    \log\log y+O(1),
\]
and
\[
    \sum_{p\le y}
    \frac{
        \cos((t-\ell/\log y)\log p)
    }{
        p^{1+\beta}
    }
    \ge
    \log\log y+O(1).
\]
Consequently,
\[
\mathbb E
|F_y(1/2+\beta+it)|^2
\left|
F_y\left(
    \frac12+\frac{i\ell}{\log y}
\right)
\right|^{2(k-1)}
\gg_k
(\log y)^{k^2}.
\]
The intersection of $[-1/2,1/2]$ with
$|t-\ell/\log y|\le1/\log y$ has length at least $1/\log y$,
including when the centre is an endpoint. Hence
\begin{equation}\label{eq:first-integral-lower}
    \mathcal I_\ell
    \gg_k
    (\log y)^{k^2-1}.
\end{equation}

We next consider the error integral
\[
\mathcal J_\ell
:=
\int_{-\infty}^{\infty}
\frac{
\mathbb E
|F_y(1/2+\beta/2+it)|^2
\left|
F_y\left(
    \frac12+\frac{i\ell}{\log y}
\right)
\right|^{2(k-1)}
}{
1/4+t^2
}\,dt.
\]
Lemma~\ref{lem:euler-correlation} gives
\begin{align*}
&\mathbb E|F_y(1/2+\beta/2+it)|^2
\left|F_y\left(\frac12+\frac{i\ell}{\log y}\right)\right|^{2(k-1)}\\
&\qquad\ll_k (\log y)^{(k-1)^2+1}
\exp\left(
2(k-1)\sum_{p\le y}
\frac{\cos((t-\ell/\log y)\log p)}{p^{1+\beta/2}}
\right).
\end{align*}
Since $\beta=1/\log y$, replacing $p^{-1-\beta/2}$ by $p^{-1}$ in an upper bound costs only $O_k(1)$ in the exponent, because
\[
\sum_{p\le y}\frac{1-p^{-\beta/2}}p
\le \frac{\beta}{2}\sum_{p\le y}\frac{\log p}{p}\ll 1.
\]
Moreover, $|\ell|/\log y\le1/2$, so $(1/4+t^2)^{-1}$ is bounded by a constant multiple of $(1+(t-\ell/\log y)^2)^{-1}$. With $u=t-\ell/\log y$, we split the integral into the ranges $|u|\le1/\log y$, $1/\log y<|u|\le10$, and $|u|>10$. Using \eqref{eq:prime-harmonic-shift}, we obtain
\begin{align*}
\mathcal J_\ell
&\ll_k (\log y)^{(k-1)^2+1}
\left(\frac{(\log y)^{2(k-1)}}{\log y}
+\int_{1/\log y}^{10}u^{-2(k-1)}\,du
+\int_{10}^{\infty}\frac{(\log u)^{2(k-1)}}{u^2}\,du\right)\\
&\ll_k(\log y)^{k^2-1},
\end{align*}
uniformly in $\ell$.

Using \eqref{eq:conditioned-final},
\eqref{eq:first-integral-lower}, and summing over the
$2\lfloor(\log y)/2\rfloor+1\asymp\log y$ values of $\ell$, we get
\begin{align}
&\sum_{|\ell|\le(\log y)/2}
\mathbb E
\left|
    \sum_{n\le x}f(n)
\right|^2
\left|
F_y\left(
    \frac12+\frac{i\ell}{\log y}
\right)
\right|^{2(k-1)}
\nonumber\\
&\quad\gg_k
\frac{x}{\log y}
(\log y)
(\log y)^{k^2-1}
\left(
1-O_k(x^{-\beta/4})
\right).
\end{align}
Since $\log x\ge C_0\log y$, we have $x^{-\beta/4}\le e^{-C_0/4}$. Thus the bracket is at least $1/2$ for
all $C_0$ above a threshold depending only on $k$. In particular the
lower-bound constant is uniform in such $C_0$.

Thus
\[
\sum_{|\ell|\le(\log y)/2}
\mathbb E
\left|
    \sum_{n\le x}f(n)
\right|^2
\left|
F_y\left(
    \frac12+\frac{i\ell}{\log y}
\right)
\right|^{2(k-1)}
\gg_k
x(\log y)^{k^2-1}.
\]
Combining this with \eqref{eq:exp-to-Euler} proves
\eqref{eq:random-main-complete}.
\end{proof}

\subsubsection{The truncation error}

Recall
\[
    I_m=(y_{m-1},y_m]
\]
and define
\[
    \mathcal P_m
    :=
    \{p\le y:p\notin I_m\}.
\]
We also put
\[
    \operatorname{Err}_{m,\ell}(f)
    :=
    \exp(2(k-1)\Re D_{m,\ell}(f))
    -
    R_{m,\ell}(f).
\]

The following variant of \cite[Lemma~6]{Szabo25} is uniform without
the factor $\log x/\log y$ in that estimate. This matters in the long
range, where $y=(T/x)^{1/C_0}$.

\begin{lemma}\label{lem:weighted-block-moment}
Let $j\ge0$ be an integer and put
\[
    B_m
    :=
    \sum_{p\in I_m}
    \left(
        \frac2p+\frac3{p^2}
    \right).
\]
Then uniformly in $t\in\mathbb R$,
\begin{align}
&\mathbb E
\left|
    \sum_{n\le x}f(n)
\right|^2
\exp\left(
    2(k-1)\Re
    \sum_{p\in\mathcal P_m}
    \left(
        \frac{f(p)}{p^{1/2+it}}
        +
        \frac{f(p)^2}{2p^{1+2it}}
    \right)
\right)
\nonumber\\
&\qquad\qquad\times
|D_{m}(f;t)|^{2j}
\nonumber\\
&\quad\ll_k
x(\log y)^{k^2-2}
j!\,B_m^{j}
\prod_{p\in I_m}
\left(1-\frac1p\right)^{-2},
\label{eq:weighted-block-moment}
\end{align}
where
\[
    D_m(f;t)
    :=
    \sum_{p\in I_m}
    \left(
        \frac{f(p)}{p^{1/2+it}}
        +
        \frac{f(p)^2}{2p^{1+2it}}
    \right).
\]
\end{lemma}

\begin{proof}
Put $\alpha=k-1\ge1$, $\mathcal P=\mathcal P_m$, and write
$\langle\mathcal P\rangle$ for the positive integers supported on
$\mathcal P$, including $1$. Set
\[
  V(t)=\exp\!\left(2\alpha\Re\sum_{p\in\mathcal P}
    \left(\frac{f(p)}{p^{1/2+it}}+
    \frac{f(p)^2}{2p^{1+2it}}\right)\right),\qquad
  A(X)=\sum_{\substack{d\le X\\p\mid d\Rightarrow p\notin I_m}}f(d).
\]
Let $\mathcal F$ be the sigma-algebra generated by $f(p)$ for all primes
$p\le x$ outside $I_m$. Both $V(t)$ and $A(X)$ for $X\le x$ are
$\mathcal F$-measurable. Each integer factors uniquely into its part
supported on $I_m$ and its part free of primes in $I_m$. Therefore
\[
 \sum_{n\le x} f(n)
    =\sum_{\substack{r\le x\\p\mid r\Rightarrow p\in I_m}}f(r)A(x/r).
\]
Conditional on $\mathcal F$, the coefficients $A(x/r)$ are fixed
and only the independent phases in $I_m$ remain random. Applying
\cite[Lemma~4]{Szabo25} in these phases and taking expectation over
$\mathcal F$ bounds the left side of \eqref{eq:weighted-block-moment}
by an absolute constant times
\begin{equation}\label{eq:conditional-block}
 j! B_m^j \sum_{\substack{r\le x\\p\mid r\Rightarrow p\in I_m}}
  \tau(r)\,\mathbb E[V(t)|A(x/r)|^2].
\end{equation}
Here $\tau(r)=\sum_{d\mid r}1$, because every divisor of $r$ is supported
on the block, as required by the divisor weight in the cited lemma.
The two terms at a prime in
that inequality give $2/p+6/(4p^2)\le 2/p+3/p^2$. For $j=0$ the
same estimate follows directly by conditional orthogonality.

We claim, uniformly in $t$ and $1\le X\le x$, that
\begin{equation}\label{eq:smooth-weighted-length}
  \mathbb E[V(t)|A(X)|^2]\ll_k X(\log y)^{k^2-2}.
\end{equation}
Let $d_\alpha$ denote the coefficients of $\zeta(s)^\alpha$. The
coefficients of $\exp(\alpha(z+z^2/2))$ are non-negative and bounded
termwise by those of $(1-z)^{-\alpha}$. Expanding the squared modulus
in $V(t)|A(X)|^2$ and taking absolute values of the deterministic
coefficients, Steinhaus orthogonality gives the following upper bound
after dropping the restrictions on prime factors of $d,e$ from $A(X)$:
\[
 \sum_{u,v\in\langle\mathcal P\rangle}
  \frac{d_\alpha(u)d_\alpha(v)}{\sqrt{uv}}
  \sum_{d,e\le X}\mathbf1_{du=ev}.
\]
The expansions converge absolutely since $\mathcal P$ is finite. Write
$u=ga$, $v=gb$, $(a,b)=1$. Then $d=hb$, $e=ha$, so the inner count is
at most $X/\max(a,b)$. Also $d_\alpha(ga)\le d_\alpha(g)d_\alpha(a)$,
because $d_\alpha(p^j)/d_\alpha(p^{j-1})=1+(\alpha-1)/j$ decreases
with $j$. Dropping the coprimality condition and enlarging the prime
support to all $p\le y$, we obtain
\begin{align*}
 \mathbb E[V(t)|A(X)|^2]
 &\le X\left(\sum_{P^+(g)\le y}\frac{d_\alpha(g)^2}{g}\right)
 \sum_{P^+(ab)\le y}
 \frac{d_\alpha(a)d_\alpha(b)}{\sqrt{ab}\max(a,b)}.
\end{align*}
The first factor in parentheses is $\ll_k(\log y)^{\alpha^2}$, by the
local expansion $1+\alpha^2/p+O_k(p^{-2})$. For the double sum use
\[
 \frac{\sqrt{ab}}{\max(a,b)}=e^{-|\log(a/b)|/2}
  =\frac{1}{2\pi}\int_{\mathbb R}
        \frac{(a/b)^{iu}}{1/4+u^2}\,du.
\]
Thus the double sum equals
\[
 \frac{1}{2\pi}\int_{\mathbb R}
   \left|\prod_{p\le y}(1-p^{-1-iu})^{-1}\right|^{2\alpha}
                  \frac{du}{1/4+u^2}.
\]
The logarithm of the absolute-value factor is
$2\alpha\sum_{p\le y}\cos(u\log p)/p+O_k(1)$.
By \eqref{eq:prime-harmonic-shift}, the integral is
\[
 \ll_k (\log y)^{2\alpha-1}
 +\int_{1/\log y}^{10}u^{-2\alpha}\,du
 +\int_{10}^\infty\frac{(\log u)^{2\alpha}}{u^2}\,du
 \ll_k(\log y)^{2\alpha-1}.
\]
Since $\alpha^2+2\alpha-1=k^2-2$, this proves
\eqref{eq:smooth-weighted-length}.

In \eqref{eq:conditional-block} take $X=x/r$. As $r$ is supported
entirely on the selected block, the remaining sum satisfies
\[
 \sum_{\substack{r\le x\\p\mid r\Rightarrow p\in I_m}}
    \frac{\tau(r)}{r}
 \le \prod_{p\in I_m}\sum_{h\ge0}\frac{h+1}{p^h}
 = \prod_{p\in I_m}(1-1/p)^{-2}.
\]
Combining the last two estimates with \eqref{eq:conditional-block} gives
\eqref{eq:weighted-block-moment}.
\end{proof}

\begin{proposition}\label{prop:random-truncation}
For every $1\le m\le M$,
\begin{align}
&\sum_{|\ell|\le(\log y)/2}
\mathbb E
\left|
    \sum_{n\le x}f(n)
\right|^2
\exp\left(
    2(k-1)
    \sum_{\substack{1\le m'\le M\\m'\ne m}}
    \Re D_{m',\ell}(f)
\right)
|\operatorname{Err}_{m,\ell}(f)|
\nonumber\\
&\qquad\ll_k
e^{-J_m}
x(\log y)^{k^2-1}.
\label{eq:random-truncation-complete}
\end{align}
\end{proposition}

\begin{proof}
Put
\[
    K_{m,\ell}(f)
    :=
    \left|
        \sum_{n\le x}f(n)
    \right|^2
    \exp\left(
        2(k-1)
        \sum_{m'\ne m}
        \Re D_{m',\ell}(f)
    \right).
\]
By Lemma~\ref{lem:weighted-block-moment},
\begin{align}
\mathbb E
K_{m,\ell}(f)
|\Re D_{m,\ell}(f)|^{2j}
&\ll_k
x(\log y)^{k^2-2}
j! B_m^{j}
\prod_{p\in I_m}
\left(1-\frac1p\right)^{-2}.
\label{eq:even-block}
\end{align}
By Cauchy--Schwarz,
\begin{align}
\mathbb E
K_{m,\ell}(f)
|\Re D_{m,\ell}(f)|^{2j+1}
&\le
\left(
\mathbb E K_{m,\ell}
|\Re D_{m,\ell}|^{2j}
\right)^{1/2}
\nonumber\\
&\quad\times
\left(
\mathbb E K_{m,\ell}
|\Re D_{m,\ell}|^{2j+2}
\right)^{1/2},
\end{align}
so for odd moments we obtain the corresponding bound with
$(j+1)!B_m^{j+1/2}$ in place of $j!B_m^j$.

Now
\[
\operatorname{Err}_{m,\ell}(f)
=
\sum_{\substack{j_1,j_2\ge0\\
                 \max\{j_1,j_2\}>J_m}}
\frac{(k-1)^{j_1+j_2}}{j_1!j_2!}
(\Re D_{m,\ell}(f))^{j_1+j_2}.
\]
Thus
\begin{align}
\mathbb E
K_{m,\ell}
|\operatorname{Err}_{m,\ell}|
&\ll_k
x(\log y)^{k^2-2}
\prod_{p\in I_m}
\left(1-\frac1p\right)^{-2}
\nonumber\\
&\quad\times
\sum_{\substack{j_1,j_2\ge0\\
                 \max\{j_1,j_2\}>J_m}}
\frac{(k-1)^{j_1+j_2}}{j_1!j_2!}
\left\lceil\frac{j_1+j_2}{2}\right\rceil!
B_m^{(j_1+j_2)/2}.
\label{eq:tail-sum}
\end{align}

We have
\[
    B_1\ll\log\log y,
    \qquad
    B_m\ll1\quad(m\ge2),
\]
whereas
\[
    J_1=\lceil(\log\log y)^{3/2}\rceil,
    \qquad
    J_m\ge J_M\gg_k1.
\]
By first taking $C_0$ large and then $L_0(k)$ large, we have
\[
    10^4(k-1)^2B_m\le J_m.
\]
By symmetry we may suppose $j_1\ge j_2$. Then
\[
    \left\lceil\frac{j_1+j_2}{2}\right\rceil!
    \le
    j_1^{(j_1+j_2)/2}.
\]
Hence the tail in \eqref{eq:tail-sum} is
\begin{align*}
&\ll
\sum_{j_1>J_m}
\frac{(k-1)^{j_1}}{j_1!}
(B_mj_1)^{j_1/2}
\sum_{j_2\ge0}
\frac{(k-1)^{j_2}}{j_2!}
(B_mj_1)^{j_2/2}\\
&=
\sum_{j_1>J_m}
\frac{(k-1)^{j_1}}{j_1!}
(B_mj_1)^{j_1/2}
\exp\left(
    (k-1)(B_mj_1)^{1/2}
\right).
\end{align*}
The parameter condition implies
\[
    (k-1)(B_mj_1)^{1/2}
    \le j_1
\]
and, using $j_1!\ge(j_1/e)^{j_1}$,
the last expression is
\[
    \ll
    \sum_{j_1>J_m}
    \left(
        \frac{
        e^2(k-1)B_m^{1/2}
        }{
        j_1^{1/2}
        }
    \right)^{j_1}
    \ll
    e^{-2J_m}.
\]
Furthermore,
\[
    \prod_{p\in I_m}
    \left(1-\frac1p\right)^{-2}
    \ll e^{J_m}.
\]
Thus
\[
\mathbb E
K_{m,\ell}
|\operatorname{Err}_{m,\ell}|
\ll_k
e^{-J_m}
x(\log y)^{k^2-2}.
\]
There are $O(\log y)$ values of $\ell$, which gives
\eqref{eq:random-truncation-complete}.
\end{proof}

\subsubsection{The mixed proxy moment}

We now prove the random-model estimate needed for the denominator.

Recall the intervals $I_n^{(m)}$, the quantities $W_m$, and the majorants $U_{m,\ell}^{(n_m)}$ from \eqref{eq:def-U}. We write $U_{m,\ell}(f)$ for the same polynomial with $D_{m,\ell}(t)$ replaced by $D_{m,\ell}(f)$.

\begin{lemma}\label{lem:good-mixed-block}
If $n_m=0$, then
\begin{align}
\mathbb E\left|
\exp\left(2(k-1)\Re D_{m,\ell_1}(f)+2\Re D_{m,\ell_2}(f)\right)
-R_{m,\ell_1}(f)U_{m,\ell_2}(f)
\right|
\ll_k e^{-J_m}.
\label{eq:good-block}
\end{align}
\end{lemma}

\begin{proof}
Expanding both truncated exponentials, the difference between
\[
R_{m,\ell_1}(f)U_{m,\ell_2}(f)
\]
and
\[
\exp\left(
    2(k-1)\Re D_{m,\ell_1}(f)
    +
    2\Re D_{m,\ell_2}(f)
\right)
\]
is bounded by
\begin{align*}
\sum_{\substack{j_1,j_2,j_3,j_4\ge0\\
\max j_r>J_m}}
\frac{(k-1)^{j_1+j_2}}
{j_1!j_2!j_3!j_4!}
|D_{m,\ell_1}|^{j_1+j_2}
|D_{m,\ell_2}|^{j_3+j_4}.
\end{align*}
Write $s=j_1+j_2+j_3+j_4>J_m$ and $r=\lceil s/2\rceil$.
For independent Steinhaus phases $Z_p$, the multinomial expansion gives
$\mathbb E|\sum b_p Z_p|^{2r}\le r!(\sum|b_p|^2)^r$.
Apply this separately to the $f(p)$ and the $f(p)^2$ terms, whose phases
are each independent across primes. Minkowski's inequality and
$(u+v)^2\le2u^2+2v^2$ give
$\mathbb E|D_{m,\ell}|^{2r}\le r! B_m^r$.
Monotonicity of norms and H\"older's inequality therefore bound the
expectation of the product of the two powers of $|D|$ by
$(C\sqrt{sB_m})^s$, with $C$ absolute. By the multinomial formula, the
sum of the coefficients for fixed $s$ is at most $(2k)^s/s!$.
Thus the expected tail is
\[
 \ll\sum_{s>J_m}\frac{(2k C\sqrt{s B_m})^s}{s!}
 \le \sum_{s>J_m}\left(C_k\sqrt{B_m/s}\right)^s
 \ll e^{-J_m}.
\]
Indeed, $B_1/J_1\ll (\log\log y)^{-1/2}$ and
$\max_{m\ge2} B_m/J_m\ll J_M^{-1}$. The parameter choices make
$C_k\sqrt{B_m/J_m}\le e^{-2}$ for every $m$.
\end{proof}

\begin{lemma}\label{lem:bad-mixed-block}
If $n_m\ge1$, then
\begin{equation}\label{eq:bad-block}
    \mathbb E
    R_{m,\ell_1}(f)U_{m,\ell_2}(f)
    \le
    (1+W_m)^{-2}.
\end{equation}
\end{lemma}

\begin{proof}
By Cauchy--Schwarz,
\[
\left(
\mathbb E
R_{m,\ell_1}(f)U_{m,\ell_2}(f)
\right)^2
\le
\mathbb ER_{m,\ell_1}(f)^2
\,
\mathbb EU_{m,\ell_2}(f)^2.
\]

First,
\[
R_{m,\ell_1}(f)^2
\le
\exp\left(
    4(k-1)|D_{m,\ell_1}(f)|
\right).
\]
Using
\[
    e^u\ll\cosh u
    =
    \sum_{j\ge0}\frac{u^{2j}}{(2j)!}
\]
and
\[
    \mathbb E|D_{m,\ell_1}|^{2j}
    \ll j!B_m^j,
\]
we obtain
\[
    \mathbb ER_{m,\ell_1}(f)^2
    \ll_k
    e^{O(k^2B_m)}.
\]

Write $J=J_m$, $a=a_m$, $W=W_m$ and $W_0=100kJ$.
For every real $J/(100k)\le W\le W_0$ let
$U_{\rm mid}(W)=e^{4W}|D_{m,\ell_2}/W|^a$. The prime-block bound
$\mathbb E|D_{m,\ell_2}|^{2a}\le a!B_m^a$, together with
$8W\le2a$ and $a!\le a^a$, gives
\[
  \mathbb E U_{\rm mid}(W)^2
    \le \left(\frac{e^2 a B_m}{W^2}\right)^a.
\]
Uniformly in this range, $aB_m/W^2\ll_k B_m/J$. This can be made
arbitrarily small for all blocks: for $m\ge2$ use $B_m\ll1$ and
$J\ge J_M$, and for $m=1$ use $B_1/J_1\ll(\log\log y)^{-1/2}$.
We can thus ensure
$\mathbb E U_{\rm mid}(W)^2\le e^{-2a}$ and
$\mathbb E R_{m,\ell_1}^2\le e^a$, by the estimate above.
Cauchy--Schwarz now gives
\begin{equation}\label{eq:mid-endpoint}
 \mathbb E[R_{m,\ell_1}U_{\rm mid}(W)]
 \le e^{-a/2}\le(1+W)^{-2}
 \qquad (J/(100k)\le W\le W_0),
\end{equation}
since $a\ge4W_0$ and $W_0$ is sufficiently large.

For $W>W_0$, denote the last expression in \eqref{eq:def-U} by
$U_*(W)$ and evaluate the same formula at $W_0$. With
$z=2(k-1)W_0\ge2(J+1)$ we have
$2z^J/J!\le z^{J+1}/(J+1)!\le e^z$. It follows that
$U_*(W_0)\le U_{\rm mid}(W_0)$. Moreover
\[
 U_*(W)=\left(\frac{W_0}{W}\right)^{a-2J/(k-1)}U_*(W_0).
\]
Since $a-2J/(k-1)\ge2$, equation \eqref{eq:mid-endpoint} at $W_0$
yields
\[
 \mathbb E[R_{m,\ell_1}U_*(W)]
  \le\left(\frac{W_0}{W}\right)^2(1+W_0)^{-2}
  \le(1+W)^{-2}.
\]
This proves \eqref{eq:bad-block} uniformly in the cell.
\end{proof}

\begin{proposition}\label{prop:mixed-random-complete}
Uniformly for
\[
    |\ell_1|,|\ell_2|
    \le\frac{\log y}{2},
\]
and for every fixed tuple $(n_1,\ldots,n_M)$,
\begin{equation}\label{eq:mixed-random-complete}
\mathbb E
\prod_{m=1}^M
R_{m,\ell_1}(f)U_{m,\ell_2}(f)
\ll_k
\frac{
    (\log y)^{k^2}
}{
    1+|\ell_1-\ell_2|^{2(k-1)}
}
\prod_{m=1}^M(1+W_m)^{-2}.
\end{equation}
\end{proposition}

\begin{proof}
Since the prime intervals $I_m$ are disjoint, the random variables
\[
    R_{m,\ell_1}(f)U_{m,\ell_2}(f),
    \qquad 1\le m\le M,
\]
are independent. Hence
\[
\mathbb E
\prod_{m=1}^M
R_{m,\ell_1}(f)U_{m,\ell_2}(f)
=
\prod_{m=1}^M
\mathbb E
R_{m,\ell_1}(f)U_{m,\ell_2}(f).
\]

If $n_m\ge1$, Lemma~\ref{lem:bad-mixed-block} gives the factor $(1+W_m)^{-2}$ with no block-dependent constant.

Suppose $n_m=0$. Since the exponential in Lemma~\ref{lem:good-mixed-block} has expectation at least $1$ by Jensen's inequality, that lemma gives
\[
\mathbb E R_{m,\ell_1}(f)U_{m,\ell_2}(f)
\le (1+O_k(e^{-J_m}))
\mathbb E\exp\left(2(k-1)\Re D_{m,\ell_1}(f)+2\Re D_{m,\ell_2}(f)\right).
\]
Independence over primes and a second-order local expansion give
\begin{align}
&\mathbb E\exp\left(2(k-1)\Re D_{m,\ell_1}(f)+2\Re D_{m,\ell_2}(f)\right)\nonumber\\
&\quad=\exp\Bigg(\sum_{p\in I_m}\left\{\frac{(k-1)^2+1}{p}
+\frac{2(k-1)\cos\left((\ell_1-\ell_2)\log p/\log y\right)}p
+O_k(p^{-3/2})\right\}\Bigg).
\label{eq:mixed-local}
\end{align}
The products of $(1+O_k(e^{-J_m}))$ over the good blocks are bounded, and the sum of the $O_k(p^{-3/2})$ terms over all blocks is $O_k(1)$. Notice also that
\[
(k-1)^2+1+2(k-1)\cos\theta\ge (k-2)^2\ge0.
\]
Thus adding the primes belonging to bad blocks can only enlarge the main exponential factor. Multiplying the good-block estimates and using Lemma~\ref{lem:bad-mixed-block} on the bad blocks, we obtain
\begin{align}
&\mathbb E
\prod_{m=1}^M
R_{m,\ell_1}(f)U_{m,\ell_2}(f)
\nonumber\\
&\quad\ll_k
\exp\left(
((k-1)^2+1)
\sum_{p\le y}\frac1p
+
2(k-1)
\sum_{p\le y}
\frac{
\cos\left(
    \frac{(\ell_1-\ell_2)\log p}{\log y}
\right)
}{p}
\right)
\prod_m(1+W_m)^{-2}.
\label{eq:mixed-global}
\end{align}

By Mertens' theorem,
\[
    \sum_{p\le y}\frac1p
    =
    \log\log y+O(1),
\]
while the shifted prime harmonic estimate gives
\[
\sum_{p\le y}
\frac{
\cos\left(
    \frac{(\ell_1-\ell_2)\log p}{\log y}
\right)
}{p}
\le
\log\log y
-
\log(1+|\ell_1-\ell_2|)
+
O(1).
\]
Consequently the exponential in \eqref{eq:mixed-global} is
\[
\ll_k
(\log y)^{
    (k-1)^2+1+2(k-1)
}
(1+|\ell_1-\ell_2|)^{-2(k-1)}.
\]
Since
\[
    (k-1)^2+1+2(k-1)=k^2,
\]
we conclude that
\[
\mathbb E
\prod_{m=1}^M
R_{m,\ell_1}(f)U_{m,\ell_2}(f)
\ll_k
\frac{
(\log y)^{k^2}
}{
1+|\ell_1-\ell_2|^{2(k-1)}
}
\prod_m(1+W_m)^{-2}.
\]
This proves \eqref{eq:mixed-random-complete}.
\end{proof}
	\section{The upper bounds}\label{sec:upper}
Throughout this section write $\|f\|_{p,T}=(T^{-1}\int_T^{2T}|f(t)|^p\dd t)^{1/p}$.
\subsection{Shifted moments to the right of the critical line}
For $H\ge 4$ and $|u|\le1$, put $g_H(u)=\min(H,1/|u|)$, with $g_H(0)=H$. For fixed positive $b_1,\dots,b_r$ define
\[
 B_H(\mathbf h)= H^{\sum_j b_j^2}
   \prod_{i<j}g_H(h_i-h_j)^{2b_i b_j},
 \qquad 0\le h_j\le1.
\]
Exponents equal to zero can simply be omitted.

\begin{lemma}\label{lem:upper-shifts}
Under RH, uniformly for $4\le H\le\log T$ and $\mathbf h\in[0,1]^r$,
\[
 \frac1T\int_T^{2T}\prod_{j=1}^r
 \left|\zeta\left(\tfrac12+H^{-1}+i(t+h_j)\right)\right|^{2b_j}\dd t
 \ll_{\mathbf b} B_H(\mathbf h).
\]
\end{lemma}
\begin{proof}
Write $a=1/H$, $R=(\log T)/H$, $F(s)=\prod_j\zeta(s+ih_j)^{b_j}$. At heights in $[T/2,3T]$, RH gives analytic branches on $\Re s>1/2$, chosen from the Euler product on $\Re s>1$. Each factor has modulus $|\zeta(s+ih_j)|^{b_j}$. At a zero on the boundary its limiting value is zero.

On the critical line, Curran's theorem \cite[Theorem 1.1]{Curran23} gives
\begin{equation}\label{eq:upper-curran}
 \frac1T\int_{T/2}^{3T}|F(1/2+it)|^2\dd t\ll_{\mathbf b} B_{\log T}(\mathbf h).
\end{equation}
Use finitely many dyadic intervals and $|\zeta(1+1/\log T+iu)|\ll g_{\log T}(u)$ for $|u|\le1$, by the Laurent expansion at $1$. Since $g_{\log T}(u)\le Rg_H(u)$,
\begin{equation}\label{eq:upper-ratio}
 B_{\log T}(\mathbf h)\le R^{(\sum_j b_j)^2}B_H(\mathbf h).
\end{equation}

Choose a fixed smooth non-increasing $\psi$ equal to $1$ on $(-\infty,1/4]$ and to $0$ on $[1/2,\infty)$. Define
\[
 K(w)=-\frac1w\int_{1/4}^{1/2}\psi'(v)T^{vw}\dd v.
\]
Its residue at $0$ equals $1$. Integration by parts, with the trivial bound when $|w|\log T<1$, gives, for every integer $D\ge1$,
\[
 |K(w)|\ll_D\frac{\max_{1/4\le v\le1/2}T^{v\Re w}}
 {|w|(1+|w|\log T)^D}.
\]
In particular
\begin{equation}\label{eq:upper-kernel}
 \int_{\R}|K(-a+iv)|\dd v\ll e^{-R/4},
\end{equation}
by setting $v=au$ and using $R\ge1$.

Let $F(s)=\sum_{n\ge1}c(n)n^{-s}$ for $\Re s>1$. For $s=1/2+a+it$, $T\le t\le 2T$, Mellin inversion gives
\[
 P(s):=\sum_{n\le\sqrt T}c(n)n^{-s}\psi\left(\frac{\log n}{\log T}\right)
   =\frac1{2\pi i}\int_{(1)}F(s+w)K(w)\dd w .
\]
Indeed the inverse Mellin transform of $T^{vw}/w$ is the cutoff at $T^v$. Truncate at $|\Im w|=T/4$ and move to $\Re w=-a$, crossing the pole at $0$. On the horizontal sides the heights of all zeta factors are $\asymp T$ and their moduli are polynomially bounded. Taking $D$ sufficiently large bounds these integrals and the original tails by $O_{\mathbf b}(T^{-2})$. One may first move to $-a+\delta$ and let $\delta\downarrow0$ in the finite rectangle. Minkowski's inequality, \eqref{eq:upper-curran}, and \eqref{eq:upper-kernel} yield
\[
 \|F(1/2+a+i\,\cdot)-P(1/2+a+i\,\cdot)\|_{2,T}
 \ll_{\mathbf b}e^{-R/4} B_{\log T}(\mathbf h)^{1/2}+T^{-2}
 \ll_{\mathbf b} B_H(\mathbf h)^{1/2},
\]
where the last step uses \eqref{eq:upper-ratio}.

The Dirichlet polynomial mean-value theorem \cite{MV74} gives
\[
 \|P(1/2+a+i\,\cdot)\|_{2,T}^2
 \ll\sum_{n\ge1}|c(n)|^2n^{-1-2a}.
\]
Here $c$ is multiplicative, $c(p)=\sum_jb_jp^{-ih_j}$, and $|c(p^\nu)|\le d_{\sum b_j}(p^\nu)$, where $d_b$ denotes the coefficients of $\zeta(s)^b$. The terms of degree $\nu\ge2$ in the local factors contribute $O_{\mathbf b}(1)$ to the logarithm, uniformly in $a$. Consequently the last sum is
\[
 \ll_{\mathbf b}\zeta(1+2a)^{\sum b_j^2}
 \prod_{i<j}|\zeta(1+2a+i(h_i-h_j))|^{2b_i b_j}
 \ll_{\mathbf b} B_H(\mathbf h).
\]
This proves the lemma.
\end{proof}

\subsection{Integrating the shifts for real moments}
\begin{lemma}\label{lem:upper-integral}
Suppose RH holds. For fixed real $p=2k>4$ and $4\le H\le\log T$,
\[
 \frac1T\int_T^{2T}
 \left(\int_0^1|\zeta(1/2+H^{-1}+i(t+v))|\dd v\right)^p\dd t
 \ll_p H^{(k-1)^2}.
\]
\end{lemma}
\begin{proof}
We use the shift integration of Szab\'o \cite[Sections~6--7]{Szabo24} with the logarithmic scale $H$ of Lemma~\ref{lem:upper-shifts}.
Put $b=p-3>1$ and $Z(u)=|\zeta(1/2+H^{-1}+i(t+u))|$.
Expand three copies of $\int_0^1 Z$. In the remaining integral over $u$, partition according to the six possible orders of $|u-h_i|$, $1\le i\le3$. The convexity inequality for a sum of six terms and symmetry allow us to use only $r_1\le r_2\le r_3$, where $r_i=|u-h_i|$.

Put $I_0=[0,H^{-1}]$ and $I_j=(2^{j-1}/H,2^j/H]$ for $j\ge1$. First split the $u$-integral according to $r_1\in I_j$. Weighted H\"older gives $(\sum A_j)^b\ll_b\sum(1+j)^{2b}A_j^b$ for $A_j\ge0$. For fixed $h_1,h_2,h_3$ the set has $u$-length $\ll2^j/H$, so H\"older inside this set gives an additional factor $\ll_b(2^j/H)^{b-1}$ when replacing the integrand by $Z(u)^b$. Now partition the four-variable integral by
\[
 r_1\in I_j,\qquad r_2-r_1\in I_l,\qquad r_3-r_2\in I_m.
\]
Each part has volume $\ll 2^{j+l+m}/H^3$, by fixing $u$ and the three signs of $h_i-u$. Write $M=\max(j,l,m)$. In Lemma~\ref{lem:upper-shifts} use exponents $1/2,1/2,1/2,b/2$ at $h_1,h_2,h_3,u$. Using $|h_i-h_s|\ge |r_i-r_s|$, its bound on this part is
\[
 \ll_p H^{k^2}
 2^{-\frac b2(j+\max(j,l)+M)
    -\frac12(l+m+\max(l,m))}.
\]
Combining volume and H\"older factors gives the power $H^{k^2-3-(b-1)}=H^{(k-1)^2}$. The exponent of $2$ is
\[
 \tfrac b2\big(j-\max(j,l)-M\big)+\tfrac12\min(l,m)
 \le-\tfrac{b-1}{2}M.
\]
The triple sum converges even with the factor $(1+j)^{2b}$, since $b>1$.
\end{proof}

\subsection{Sharp cutoffs on a dyadic interval}
\begin{lemma}\label{lem:upper-max}
Under RH, for fixed $k>2$, $p=2k$, $d=(k-1)^2$ and $1\le N\le2\sqrt T$,
\[
 \left\|\sup_{N\le u\le v\le 2N}
 \left|\sum_{u<n\le v}n^{-i\,\cdot}\right|\right\|_{p,T}
 \ll_k \sqrt N\,(\log(2N))^{d/p}.
\]
\end{lemma}
\begin{proof}
It suffices to consider $H=\log(2N)\ge4$. Fix $W\in C_c^\infty((1/2,3))$ equal to $1$ on $[1,2]$, and put $P_N(t)=\sum_nW(n/N)n^{-it}$. Move the Mellin integral to $\sigma=1/2+1/H$. With $\widehat W(s)=\int_0^\infty W(y)y^{s-1}\dd y$ this gives
\[
 P_N(t)=N^{1-it}\widehat W(1-it)+
 \frac1{2\pi}\int_{\R}N^{\sigma+iv}\widehat W(\sigma+iv)
                      \zeta(\sigma+i(t+v))\dd v.
\]
Here $N^\sigma\ll\sqrt N$ and $\widehat W(\sigma+iv)\ll_D(1+|v|)^{-D}$ for every fixed $D$. For $T/2\le t\le 3T$ the residue and the integral over $|v|>T/8$ are $O(T^{-2})$, by polynomial growth of zeta and increasing $D$. Divide the remaining $v$-range into unit intervals. On each, Lemma~\ref{lem:upper-integral} applies after translating $t$ and covering heights in $[T/4,4T]$ by finitely many dyadic intervals. Its hypotheses hold there for large $T$ since $N\le 2\sqrt T$. Minkowski and the summable Mellin weights give
\begin{equation}\label{eq:upper-smooth}
 \frac1T\int_{T/2}^{3T}|P_N(t)|^p\dd t\ll_k N^{p/2}H^d.
\end{equation}

To take sharp partial sums, choose a nonnegative Schwartz function $\rho$ of integral $1$ whose Fourier transform is supported in $[-c,c]$, for a sufficiently small absolute $c>0$. Such a function is obtained by squaring the modulus of a bandlimited Schwartz function. In angular frequency convention put
\[
 \eta_T=\mathbf 1_{[3T/4,9T/4]}*\big(T^{-1/2}\rho(\,\cdot\,/\sqrt T)\big).
\]
Then $\eta_T\gg1$ on $[T,2T]$, $0\le\eta_T\le1$, and $\widehat\eta_T$ is supported in $[-c/\sqrt T,c/\sqrt T]$. Outside $[T/2,3T]$,
\[
 |\eta_T(t)|\ll_D\sqrt T\left(1+\frac{T+|t|}{\sqrt T}\right)^{-D}.
\]
The packets of $\eta_TP_N$ centred at $-\log n$ are disjoint, as their frequency gaps are $\gg1/N$. For each interval $(u,v]$ its packets can therefore be selected by a single Fourier interval projection. The real-line Carleson--Hunt theorem \cite{Kovac} bounds the $L^p$ norm of the supremum of these projections by $O_p(\|\eta_TP_N\|_{L^p(\R)})$ (an interval is the difference of two half-lines). On the discarded tails $|P_N(t)|\ll N$, so the displayed decay makes their contribution $O(1)$ to the unnormalized $p$-moment for large enough $D$. Applying \eqref{eq:upper-smooth} on $[T/2,3T]$ and dividing by $T$ completes the proof.
\end{proof}

\subsection{The transition at the square root}
\begin{proof}[Proof of Theorem~\ref{thm:upper}]
Put $d=(k-1)^2$ and $p=2k$. When $2\le x\le\sqrt T$, sum Lemma~\ref{lem:upper-max} over dyadic intervals up to $x$, including $n=1$ separately. The geometric sum is $\ll_k\sqrt{x}(\log(2x))^{d/p}$.

Suppose $\sqrt T\le x\le T$. The approximate functional equation in the form recorded by Harper \cite[p.~3]{Harper23} gives, for $T\le t\le2T$,
\begin{align}
 |S_x(t)|\ll {}& |S_{\sqrt T}(t)|+1+\frac{x}{\sqrt T}\notag\\
 &+\sqrt T\left|\sum_{t/(2\pi x)<n\le t/(2\pi\sqrt T)}
                         \frac{n^{it}}{n}\right|.\label{eq:upper-flip}
\end{align}
Indeed, put $D_z(t)=\sum_{t/(2\pi z)<n\le t/\pi} n^{-1+it}$. Harper's formula gives both $|S_z|\ll\sqrt t|D_z|+1+z/\sqrt t$ and $\sqrt t|D_z|\ll|S_z|+1+z/\sqrt t$. Apply the first with $z=x$, the second with $z=\sqrt T$, and subtract the dual tails.

Write $L=T/x$. Each dyadic block of the sum in \eqref{eq:upper-flip} has $N\gg L$ and $N\ll\sqrt T$, with $n=1$ treated separately when $L\le2\pi$. Partial summation bounds its absolute value by $O(1/N)$ times the maximal sum of Lemma~\ref{lem:upper-max}, since conjugation preserves its modulus. Therefore its contribution, after multiplying by $\sqrt T$ and summing blocks in $L^p$, is
\[
 \ll_k \sqrt T\sum_{\substack{N\ {\rm dyadic}\\N\gg L}}
       N^{-1/2}(\log(2N))^{d/p}
 \ll_k \sqrt x\,(\log(2T/x))^{d/p}.
\]
The $n=1$ contribution also satisfies this bound. The error $1+x/\sqrt T$ is $O(\sqrt x)$. Finally the short-sum result gives
$\|S_{\sqrt T}\|_{p,T}\ll_k T^{1/4}(\log T)^{d/p}$.
For $x\le T^{3/4}$ use $\log(2T/x)\gg\log T$. For $x>T^{3/4}$ use $T^{1/4}(\log T)^{d/p}\ll_k \sqrt x$. Raising the resulting bound to the power $p$ proves the theorem.
\end{proof}

\subsection{Integral moments without RH}\label{sec:integer-moments}
For $k\in\mathbb N$ let
\[
 r_{k,x}(n)=\#\{(n_1,\ldots,n_k)\in\mathbb N^k:
               n_j\le x,\ n_1\cdots n_k=n\},
 \qquad D_k(x)=\sum_{n\le x^k}r_{k,x}(n)^2.
\]
Steinhaus orthogonality gives $D_k(x)=\E|\sum_{n\le x}f(n)|^{2k}$.
By \cite[Theorem~3]{HNR}, for fixed $k$ and $x\to\infty$,
\[
 D_k(x)\sim C_k x^k(\log x)^{(k-1)^2},
\]
where the leading constant is
\[
 C_k=\binom{2k-2}{k-1}\frac{\operatorname{vol}(\mathcal B_k)}{k^{k-1}}
 \prod_{p}\left[
 \left(1-\frac1p\right)^{k^2}
 \sum_{j=0}^{\infty}\binom{j+k-1}{k-1}^{\!2}p^{-j}\right].
\]
The polytope
$\mathcal B_k$ consists of nonnegative $k\times k$ matrices whose row
and column sums are all $1$. Its volume is Euclidean in its
$(k-1)^2$-dimensional affine hull, with $\operatorname{vol}(\mathcal B_1)=1$.
In particular $C_1=1$. For $k=2$ the polytope is a segment of length $2$
and the Euler product is $1/\zeta(2)$, giving $C_2=12/\pi^2$.

On expanding $S_x(t)^k=\sum_{n\le x^k}r_{k,x}(n)n^{-it}$, the
mean-value theorem \cite{MV74} on $[T,2T]$ gives
\begin{align}
 M_k(x,T)
 &= D_k(x)+O\!\left(\frac1T\sum_{n\le x^k}n\,r_{k,x}(n)^2\right)
 \notag\\
 &=D_k(x)\left(1+O\!\left(\frac{x^k}{T}\right)\right).
 \label{eq:integer-transfer}
\end{align}
This proves \eqref{eq:integer-asymptotic}. In particular it applies
when $x\to\infty$ and $x\le T^{1/k-\delta}$ for fixed $0<\delta<1/k$.

\begin{proof}[Proof of Theorem~\ref{thm:integer}]
Fix an integer $m\ge2$, put $d=(m-1)^2$, and take $W$ and $P_N$ as in
Lemma~\ref{lem:upper-max}. Expanding the $m$th power of $P_N$ gives a
Dirichlet polynomial of length at most $(3N)^m$. Its squared coefficient
sum is at most $\|W\|_\infty^{2m} D_m(3N)
\ll_m N^m(\log(2N))^d$ by \cite[Theorem~3]{HNR}. The Dirichlet polynomial
mean-value theorem \cite{MV74} therefore gives
\[
 \frac1T\int_{T/2}^{3T}|P_N(t)|^{2m}\dd t
 \ll_m (1+N^m/T)N^m(\log(2N))^d .
\]
Applying the Fourier projection argument from Lemma~\ref{lem:upper-max}
with this estimate in place of \eqref{eq:upper-smooth} gives, for
$1\le N\le2\sqrt T$,
\begin{equation}\label{eq:upper-integer-max}
 \left\|\sup_{N\le u\le v\le2N}\left|\sum_{u<n\le v}n^{-i\,\cdot}\right|
 \right\|_{2m,T}
 \ll_m (1+N^m/T)^{1/(2m)}\sqrt N\,(\log(2N))^{d/(2m)} .
\end{equation}
By \eqref{eq:integer-transfer} we also have
$M_m(x,T)\ll_m (1+x^m/T)x^m(\log(2x))^d$. This proves the upper bound
for $x\le T^{1/m}$.

For $\sqrt T\le x\le T$, divide \eqref{eq:upper-integer-max} by $N$
and use partial summation on each dual block of \eqref{eq:upper-flip}.
When $m=2$, all these blocks satisfy $N^2\ll T$. The summation in the
proof of Theorem~\ref{thm:upper} then gives
$\|S_x\|_{4,T}\ll\sqrt{x}(\log(2T/x))^{1/4}$, using the short bound
at $x=\sqrt T$ for the initial sum. Thus the fourth-moment upper bound
holds on the full range.

For $m\ge3$, put $c=d/(2m)$. Each weighted block has norm at most
$O_m((N^{-1/2}+T^{-1/(2m)})(\log(2N))^c)$. The first term sums
geometrically from $N\gg T/x$. The second sums over $O(\log T)$
blocks. The preceding direct mean-value bound also gives
$\|S_{\sqrt T}\|_{2m,T}\ll_m T^{1/2-1/(2m)}(\log T)^c$.
Including the error in \eqref{eq:upper-flip}, we obtain
\[
 \|S_x\|_{2m,T}
  \ll_m \sqrt{x}(\log(2T/x))^c
       +T^{1/2-1/(2m)}(\log T)^{1+c}.
\]
For $x\ge T^{1-1/m+\varepsilon}$ the second term is
$O_{m,\varepsilon}(\sqrt{x})$. The lower bounds in the statement
follow from Theorem~\ref{thm:main}.
\end{proof}

\section*{Acknowledgements}
	The authors would like to thank Professor Radziwiłł for drawing their attention to \cite{HNR}. The first author is supported by  the National
	Natural Science Foundation of China (Grant No. 	1240011770). The second author is supported by   the National
	Natural Science Foundation of China (Grant No.1250012812). The third author is supported by the Fundamental Research Funds for the Central Universities (Grant No. 531118010622), the National
	Natural Science Foundation of China (Grant No. 1240011979) and the Hunan Provincial Natural Science Foundation of China (Grant No. 2024JJ6120).

\end{document}